\def\fracb#1#2{\frac{{\textstyle #1}}{{\textstyle #2}}}
\def\ra{\rightarrow}
\def\ms{\mapsto}

\def\Q{\hbox{\bf Q}}

\newtheorem{theo}{Th\'eor\`eme}

\newtheorem{cor}{Corollaire}
\def\bigpar{\par\bigskip}
\def\medpar{\par\medskip}

\def\thebibliography#1{\section*{Bibliographie\markboth
 {REFERENCES}{REFERENCES}}\list
 {[\arabic{enumi}]}{\settowidth\labelwidth{[#1]}\leftmargin\labelwidth
 \advance\leftmargin\labelsep
 \usecounter{enumi}}
 \def\newblock{\hskip .11em plus .33em minus -.07em}
 \sloppy
 \sfcode`\.=1000\relax}

\def\@begintheorem#1#2{\global\advance\@listdepth -1\relax
\list{}{\leftmargin 0pt\labelwidth -\labelsep}\item[{\sc #1\ #2}]\it}
\def\@endtheorem{\endlist\global\advance\@listdepth 1\relax}

\def\section{\@startsection {section}{1}{\z@}{3.5ex plus 1ex minus
.2ex}{2.3ex plus .2ex}{\large\bf}}
\def\subsection{\@startsection{subsection}{2}{\z@}{3.25ex plus 1ex minus
.2ex}{1.5ex plus .2ex}{\normalsize\bf}}
\def\subsubsection{\@startsection{subsubsection}{3}{\z@}{3.25ex plus 1ex
minus .2ex}{1.5ex plus .2ex}{\normalsize\bf}}

 \def\tableofcontents{\section*{Table des mati\`eres}
 \markboth{CONTENTS}{CONTENTS}
 \@starttoc{toc}}

\documentstyle[A4,11pt]{article}
\begin{document}
G\'eom\'etrie alg\'ebrique/ {\it Algebraic geometry.}
\medpar
{\large Rang de courbes elliptiques d'invariant donn\'e.}

Jean-Fran\c cois Mestre.
\medpar

{\bf R\'esum\'e.-} Nous montrons qu'il existe une infinit\'e de courbes
elliptiques d\'efinies sur $\Q$, d'invariant modulaire donn\'e, et de
rang $\geq 2$.
De plus, il existe une infinit\'e de courbes d\'efinies sur $\Q$, d'invariant
nul (resp. \'egal \`a $1728$), et de rang $\geq 6$ (resp. $\geq 4$).

\medpar
{\large On the rank of elliptic curves with given modular invariant.}

{\bf Abstract.-} We prove that there exist infinitely many elliptic curves over
$\Q$
with given modular invariant, and rank $\geq 2$. Furthermore, there exist
infinitely many elliptic curves over $\Q$ with invariant equal to $0$
(resp. $1728$),
and rank $\geq 6$ (resp. $\geq 4$).
\vspace{5ex}

Soit $k$ un corps de caract\'eristique nulle, et $t$ une
ind\'etermin\'ee.
Nous prouvons ici les th\'eor\`emes suivants:
\begin{theo} Soit $j$ un \'el\'ement de $k$.
Il existe une courbe elliptique d\'efinie sur $k(t)$,
d'invariant modulaire $j$, qui n'est
pas $k(t)$-isomorphe \`a une courbe elliptique d\'efinie sur $k$,
et qui poss\`ede deux points rationnels sur $k(t)$
lin\'eairement ind\'ependants.
\end{theo}

\begin{theo}
Il existe une courbe
elliptique d\'efinie sur $k(t)$,  dont l'invariant modulaire est
\'egal \`a $1728$ (resp. $0$),
qui n'est pas $k(t)$-isomorphe \`a
une courbe d\'efinie sur $k$,
et qui poss\`ede $4$ (resp. $6$) points rationnels sur $k(t)$
lin\'eairement ind\'ependants.
\end{theo}

On en d\'eduit par sp\'ecialisation les corollaires suivants:
\begin{cor} Soit $j$ un \'el\'ement de $\Q$. Il existe une
infinit\'e de courbes elliptiques d\'efinies sur $\Q$, non
deux \`a deux $\Q$-isomorphes, d'invariant modulaire $j$, dont
le rang du groupe de Mordell-Weil est $\geq 2$.
\end{cor}
 \begin{cor}
 Il existe une infinit\'e de courbes elliptiques d\'efinies
 sur $\Q$, d'invariant modulaire \'egal \`a $1728$ (resp. $0$),
non deux \`a deux $\Q$-isomorphes, dont le rang du groupe de
Mordell-Weil est $\ge 4$ (resp. $\ge 6$).
\end{cor}
\section{D\'emonstration du th\'eor\`eme $1$}
\begin{theo}
Soient $k$ un corps de caract\'eristique nulle, et $E$ et $E'$
deux courbes elliptiques d\'efinies sur $k$. On suppose que les
invariants modulaires $j(E)$ et $j(E')$ ne sont pas simultan\'ement
\'egaux \`a $0$ ou \`a $1728$.
Il existe alors une courbe $C$, rev\^etement quadratique de la droite
projective, d\'efinie sur $k$, et deux morphismes ind\'ependants
$p:\;\;C\ra E$ et $p':\;\;C\ra E'$
d\'efinis sur $k$.
\end{theo}

(On rappelle que deux morphismes $p:\;\;C\ra E$ et $p':\;\;C\ra E'$
sont dits ind\'ependants si les images r\'eciproques par $p^*$ et
$p'^*$ des formes de premi\`ere esp\`ece de $E$ et $E'$ sont
lin\'eairement ind\'ependantes.)
\medpar
Soient $y^2=x^3+ax+b$ une \'equation de $E$ et $y^2=x^3+a'x+b'$
une \'equation de $E'$.
L'hypoth\`ese sur $j(E)$ et $j(E')$ implique que
$a=0\rightarrow a'\neq 0$ et
$b=0\rightarrow b'\neq 0$.

Posons $f(x)=x^3+ax+b$ et $g(x)=x^3+a'x+b'$.
Si $u$ est une ind\'etermin\'ee,
l'\'equation (en $x$) $$u^6f(x)=g(u^2x)$$ a pour solution
$x=\phi(u)$, avec
$\phi(u)=-\fracb{b'-u^6b}{u^2(a'-u^4a)}.$

Soit $C$ la courbe d'\'equation
$Y^2=f(\phi(X)).$
Soient
$\rho:\;\;\;C\rightarrow E$ et
$\rho':\;\;\;C\rightarrow E'$ les morphismes donn\'es par
$\rho(X,Y)=(x=\phi(X),y=Y)$ et $\rho'(X,Y)=(x=X^2\phi(X),y=X^3 Y)$.
Si $\omega=\rho^*(dx/y)$ et $\omega'=\rho'^*(dx/y)$, on a

$$\omega/\omega'=
\frac{{ { 3 a { X^{4}} b'} { -  2 { X^{6}} b a'} {- b'a' }}}{{ { X^{3}}
{( { { { X^{6}} b a} { -  3 { X^{2}} b a'}+ { 2 a b'}} )}}},$$
fraction rationnelle en $X$ non constante. Par suite, $\omega$ et $\omega'$
sont ind\'ependantes dans l'espace des formes diff\'erentielles de
premi\`ere esp\`ece de $C$,  d'o\`u le th\'eor\`eme.

\medpar
{\sc Remarque.} Le calcul montre que le genre de $C$ est
$\leq 10$. Plus pr\'ecis\'ement, si l'invariant modulaire $j(E)$ de $E$
n'est pas \'egal \`a $j(E')$, et si $j(E)$ et $j(E')$ sont distincts
de $0$ et $1728$, le genre de $C$ est \'egal \`a $10$. Si $j(E)=j(E')$,
et distinct de $0$ et $1728$, le genre
de $C$ est \'egal \`a $6$. Si $j(E)=1728$, et $j(E')\neq 0$, le
genre de $C$ vaut $7$. Si $j(E)=0$, et $j(E')\neq 1728$, le genre
de $C$ vaut $8$. Enfin, si $j(E)=0$ et $j(E')=1728$, le genre de $C$
vaut $5$.

\begin{theo}
Soient $k$ un corps de caract\'eristique nulle,
et $j$ un \'el\'ement de $k$.
Il existe une courbe $C$ d\'efinie sur $k$, rev\^etement quadratique de
la droite projective, une courbe elliptique $E$ d\'efinie sur $k$
d'invariant $j$, et deux morphismes ind\'ependants
$p$ et $p'$ de $C$ dans $E$ d\'efinis sur $k$.
\end{theo}

Si $j\in k$, $j\neq 0,1728$, et si $a=b=\fracb{27j}{4(j-1728)}$, la courbe
elliptique $E$, d\'efinie sur $k$, d'\'equation $y^2=x^3+ax+b$ a comme
invariant modulaire $j$.
Le th\'eor\`eme pr\'ec\'edent permet donc de conclure, sauf si
$j=0$ ou $j=1728$.

Or la jacobienne de la courbe de genre $2$, d\'efinie sur
$\Q$, d'\'equation $y^2=x^6+1$ est $\Q$-isog\`ene
au produit de la courbe elliptique $y^2=x^3+1$, d'invariant modulaire
\'egal \`a $0$, avec elle-m\^eme. D'o\`u le r\'esultat si $j=0$.

De m\^eme, soit $C$ la courbe de genre $2$ d'\'equation
$y^2=(t^2+1)(t^2-2)(2t^2-1).$

Les morphismes $(t,y)\ms (t^2,y)$ et $(t,y)\ms (1/t^2,y/t^3)$
d\'efinissent deux rev\^etements de $C$ sur la courbe elliptique
d'\'equation $y^2=(x+1)(x-2)(2x-1)$, dont l'invariant modulaire vaut
$1728$. Cela ach\`eve la d\'emonstration du th\'eor\`eme.
\medpar
{\sc Remarques.-}  Si $E$ est une courbe elliptique d\'efinie sur
$k$, il est parfois possible de trouver une courbe hyperelliptique
d\'efinie sur $k$, de genre $<10$, dont la jacobienne est
$k$-isog\`ene \`a $E\times E\times A$, o\`u $A$ est une vari\'et\'e
ab\'elienne convenable. Par exemple:

1) Soit $E$ une courbe elliptique d\'efinie sur $k$
d'\'equation
$y^2=x^3-ax+b$, o\`u $a$ est non nul et de la forme $\alpha^2+3\beta^2$,
$\alpha$, $\beta \in k$. La conique $x_1^2+x_1x_2+x_2^2=a$ est alors
$k$-isomorphe \`a la droite projective, d'o\`u l'existence
de deux fractions rationnelles $x_1(t)$ et $x_2(t)$ telles
que la fraction rationnelle $f(t)=x_1^3-ax_1+b$  soit \'egale \`a la
fraction rationnelle $x_2^3-ax_2+b$. On en d\'eduit $2$ applications
rationnelles
$(t,y)\ms (x_1(t),y)$ et $(t,y)\ms (x_2(t),y)$
de la courbe $C$ d'\'equation $y^2=f(t)$ sur $E$. Les fractions
rationnelles $x'_1$ et $x'_2$ n'\'etant pas proportionnelles,
et la courbe $C$ \'etant de genre $3$,
on en d\'eduit que  la jacobienne de la courbe
$C$  est $k$-isog\`ene \`a $E\times E\times E_1$,
o\`u $E_1$ est une courbe elliptique d\'efinie sur $k$.
\medpar
2) Soient $E_1$ et $E_2$ deux courbes elliptiques, d\'efinies sur
$k$, dont les points d'ordre $2$ appartiennent \`a $k$.
Si $y^2=(x-a)(x-b)(x-c)$ (resp. $y^2=(x-a')(x-b')(x-c')$) est une
\'equation de $E_1$ (resp. $E_2$),  quitte \`a permuter les r\^oles
de $a,b,c$,
on peut trouver une application
affine $x\ms h(x)=\alpha x+\beta$ telle que $h(a)=a'$, $h(b)=b'$, et
$h(c)\neq c'.$ La jacobienne de la courbe de genre $2$ d'\'equations
$$y^2=(x-a)(x-b)(x-c),\;\;\;z^2=\alpha (x-a)(x-b)(x-h^{-1}(c'))$$
est alors isog\`ene \`a $E_1\times E_2$.
\bigpar
Le th\'eor\`eme 1 d\'ecoule ais\'ement du th\'eor\`eme pr\'ec\'edent.  En
effet,
si $j\in k$,
d'apr\`es le th\'eor\`eme pr\'ec\'edent,
il existe une courbe
$C$, d\'efinie sur $k$, rev\^etement quadratique de la droite projective,
une courbe elliptique $E$ d\'efinie sur $k$ d'invariant $j$, et
deux morphismes
ind\'ependants $p_1$ et $p_2$ de $C$ sur $E$.
Soit $w$ l'involution
hyperelliptique de $C$; les morphismes $p_1\circ w+p_1$ et
$p_2\circ w+p_2$ de $C$ dans $E$ sont  constants, car $w$ agit sur
la jacobienne de $C$ comme  $-1$. Par suite, les morphismes
$p'_1=p_1\circ w-p_1$ et $p'_2=p_2\circ w-p_2$ sont ind\'ependants;
si $y^2=f(t)$ est une \'equation de $C$,
et si $E_w$ est la courbe obtenue \`a partir de $E$ par
torsion par $\sqrt{f(t)}$, les points
$P_1=p'_1(t,\sqrt{f(t)})$ et $P_2=p'_2(t,\sqrt{f(t)})$ sont des
points ind\'ependants de $E_w$, rationnels sur $k(t)$.
D'o\`u le th\'eor\`eme $1$.

\section{D\'emonstration du th\'eor\`eme 2}
\subsection{Le cas des courbes d'invariant $j=1728$}

Soit $p(x)=x^4+a_2x^2+a_1x+a_0$ un \'el\'ement de $k[x]$,
dont les racines $x_i$, $1\leq i\leq 4$, appartiennent \`a $k$, et sont de
somme nulle.
La courbe $E$
d'\'equation $x^4+a_2y^2+a_1y+a_0=0$ poss\`ede $4$ points
$k$-rationnels naturels, \`a savoir les points $P_i=(x_i,x_i)$.
Si $a_0=-u^4$, o\`u $u\in k$, $E$ poss\`ede un nouveau point
$k$-rationnel, \`a savoir le point $O=(-u,0)$. Si $a_2(a_1^2-4a_0a_2)\neq 0$,
la courbe $E$ est de genre $1$, et d'invariant modulaire
\'egal \`a $1728$.

Or l'\'equation $a_0=-u^4$ s'\'ecrit $x_1x_2x_3(x_1+x_2+x_3)=u^4.$

Comme me l'a indiqu\'e J.-P. Serre, cette \'equation
a \'et\'e \'etudi\'ee par Euler ([1], p. $660$),
qui a  exhib\'e plusieurs courbes unicursales
trac\'ees sur $S$, par exemple la courbe
$$u=1,\;\;\;x_1=t\fracb{2t^2-1}{2t^2+1},\;\;\;
x_2=\fracb{2t^2-1}{2t(2t^2+1)},\;\;\;
x_3=\fracb{4t}{2t^2-1}.$$

Soit donc $x_4=-x_1-x_2-x_3$,  o\`u les $x_i$ sont donn\'es par les
formules ci-dessus,
et soit $p=\prod (x-x_i)=x^4+a_2x^2+a_1x+a_0.$
La courbe $E$, d\'efinie sur $k(t)$, d'\'equation
$x^4+a_2y^2+a_1y+a_0$ est de genre $1$; elle est  $k(t)$-isomorphe \`a la
courbe elliptique d'\'equation $y^2=x^3+a_2(a_1^2-4a_0a_2)x$.

On v\'erifie que $a_2(a_1^2-4a_0a_2)$
n'est pas une puissance quatri\`eme dans
$k(t)$; par suite,  $E$ n'est pas $k(t)$-isomorphe \`a une courbe
d\'efinie sur $k$.
\medpar

Pour prouver que les $4$ points $P_i$ sont ind\'ependants,
le point $O$ \'etant choisi comme origine, et d\'emontrer ainsi l'assertion
du th\'eor\`eme $2$ relative aux courbes d'invariant $1728$, il suffit
de v\'erifier que, pour une valeur de $t$, les sp\'ecialisations
des points $P_i$ sont des points ind\'ependants.

Or, pour $t=1$, le calcul, \`a l'aide du logiciel gp, montre que
le d\'eterminant de la matrice des hauteurs des sp\'ecialisations
des points $P_i$ est \'egal \`a $603.61237\ldots$, et est donc non nul.

\subsection{Le cas des courbes d'invariant $0$}
Soit $p\in k[X]$ un polyn\^ome unitaire de degr\'e $6$. Il existe alors
un unique polyn\^ome unitaire $g\in k[X]$, de degr\'e 2, tel que
le polyn\^ome $r=p-g^3$ soit de degr\'e $\leq 3$.

Supposons que  les racines
$x_1,\ldots,x_6$ de $p$ soient dans $k$.  La courbe $E$ d'\'equation
$r(x)+y^3=0$ contient les $6$ points $k$-rationnels $P_i=(r(x_i),g(x_i))$,
$1\leq i\leq 6$.

De plus, si le discriminant de $r$ est non nul, la courbe $E$ est de
genre $1$ et d'invariant modulaire \'egal \`a $0$.

Si le coefficient de degr\'e $3$ de $r$ est le cube d'un \'el\'ement de $k$,
l'un des points \`a l'infini de $E$ est $k$-rationnel,
et on peut le choisir comme origine $O$ de la courbe elliptique $E$.
Nous allons montrer que, si les $x_i$ sont convenablement choisis, les
points $P_i$ sont alors ind\'ependants.

Sans nuire \`a la g\'en\'eralit\'e du probl\`eme, on peut supposer que
la somme des racines $x_i$ de $p$ est nulle.  On peut donc \'ecrire $p$
sous la forme
$p(x)=x^6+a_4x^4+a_3x^3+a_2x^2+a_1x+a_0.$
On a alors $$g(x)=x^2+a_4/3,\;\;\;r(x)=a_3x^3+(a_2-a_4^2/3)x^2+a_1x-a_4^3/27.$$

Le coefficient $a_3$ du polyn\^ome $p$ est homog\`ene de degr\'e $3$ en
les racines $x_i$ de $p$.
L'hypersurface  cubique (en les variables $u$ et $x_i$, $1\leq i\leq 5$)
d'\'equation $u^3=a_3$ poss\`ede des sous-vari\'et\'es lin\'eaires
$k$-rationnelles naturelles,  par exemple $u=0, x_1=x_2=x_3=-x_4=-x_5$.

Par des manipulations classiques, cela permet d'obtenir des courbes
unicursales trac\'ees sur cette hypersurface. On trouve par exemple
$$\begin{array}{ll}
x_1=-126(35t-19)(14t-13)(t+1),&
x_2=63(-980t^3+3549 t  - 3084 t + 1135),\\
x_3=126(35 t - 19) (14 t - 13) (t + 1),&
x_4=63(1127 t^3- 3108 t^2+ 3525 t- 988),\\
x_5=- 113876 t^3+ 265629 t^2- 259980 t + 69103,&
x_6=104615 t^3 - 293412 t^2+ 232197 t - 78364.
\end{array}$$

On obtient ainsi, par la m\'ethode d\'ecrite ci-dessous, une courbe
elliptique $E$, d\'efinie sur $k(t)$, munie de $6$ points $k(t)$-rationnels.
Cette courbe est $k(t)$-isomorphe \`a la courbe $y^2=x^3-16D$, o\`u
$D$ est le discriminant du polyn\^ome $r$.

On v\'erifie que $D$ est
un polyn\^ome irr\'eductible sur $k(t)$, et n'est donc pas une puissance
sixi\`eme. Par suite, $E$ n'est pas $k(t)$-isomorphe \`a une courbe
d\'efinie sur $k$.

Pour prouver que les points $P_i$ sont ind\'ependants,
le point $O$ \'etant choisi comme origine, il suffit de le montrer
pour une valeur convenable de $t$. Or, pour $t=1$, le d\'eterminant
de la matrice des hauteurs normalis\'ees des points $P_i$ vaut
$38462030713.186929\ldots$, et est donc non nul.
\medpar
{\sc R\'ef\'erence bibliographique}
\medpar
[1] {\sc L. Dickson}, {\it History of the theory of numbers}, vol. 2,
Chelsea $1971$.

\begin{flushright}
UFR de Math\'ematiques,
Universit\'e de Paris VII\\
2 place Jussieu,
75251 Paris Cedex 05.\end{flushright}

\end{document}